\documentclass{ifacconf}
\usepackage{graphicx}      
\usepackage{natbib}        
\usepackage{amsfonts}
\usepackage{algorithmic}
\usepackage{amsmath,amssymb} 
\usepackage{url,cite}
\usepackage{graphicx}
\usepackage{tikz,pgfplots}
\renewcommand{\qed}{\hfill$\blacksquare$}

\newcommand{\bx}{\boldsymbol{x}}

\newcommand{\bu}{\boldsymbol{u}}

\newcommand{\bQ}{\boldsymbol{Q}}

\newcommand{\bA}{\boldsymbol{A}}

\newcommand{\bB}{\boldsymbol{B}}

\newcommand{\bH}{\boldsymbol{H}}
\newcommand{\bJ}{\boldsymbol{J}}
\newcommand{\bR}{\boldsymbol{R}}

\newcommand{\st}{\mathop{\text{\normalfont s.t.}}}
\newcommand{\diag}{\mathop{\text{\normalfont diag}}}

\allowbreak
\allowdisplaybreaks

\begin{document}
\begin{frontmatter}

\title{Controllability and Observability Imply Exponential Decay of Sensitivity in Dynamic Optimization}  

\thanks[footnoteinfo]{ We acknowledge support from the Grainger Wisconsin Distinguished Graduate Fellowship.}

\author{Sungho Shin and Victor M. Zavala}

\address{Department of Chemical and Biological Engineering,\\
  University of Wisconsin-Madison, Madison, WI 53706 USA \\
  (e-mail: \{sungho.shin,victor.zavala\}@wisc.edu).}

\begin{abstract}                
We study a property of dynamic optimization (DO) problems (as those encountered in model predictive control and moving horizon estimation) that is known as exponential decay of sensitivity (EDS). This property indicates that the sensitivity of the solution at stage $i$ against a data perturbation at stage $j$ decays exponentially with $|i-j|$. {Building upon our previous results, we show that EDS holds under uniform boundedness of the Lagrangian Hessian, a uniform second order sufficiency condition (uSOSC), and a uniform linear independence constraint qualification (uLICQ). Furthermore, we prove that uSOSC and uLICQ can be obtained under uniform controllability and observability. Hence, we have that uniform controllability and observability imply EDS.} These results provide insights into how perturbations propagate along the horizon and enable the development of approximation and solution schemes. We illustrate the developments with numerical examples.
\end{abstract}


\end{frontmatter}

\section{Introduction}\label{sec:intro}

This work studies the discrete-time, dynamic optimization (DO) formulation: 
\begin{subequations}\label{eqn:mpc}
\begin{align}
  \min_{\substack{x_{0:N}\\u_{0:N-1}}}\;& \sum_{i=0}^{N-1} \ell_{i}(x_i,u_i;d_i) + \ell_N(x_N;d_N)\\
  \st\;&Tx_0 = d_{-1} \quad|\quad \lambda_{-1}\label{eqn:mpc-init}\\
       &x_{i+1} = f_i(x_i,u_i;d_i),\;i\in\mathbb{I}_{[0,N-1]} \quad|\quad\lambda_{i}.
\end{align}
\end{subequations}
Here, $N\in\mathbb{I}_{>0}$ is the horizon length; for each stage (time) $i$, $x_i\in\mathbb{R}^{n_x}$ are the states, $u_i\in\mathbb{R}^{n_u}$ are the controls, $d_i\in\mathbb{R}^{n_d}$ are the data (parameters), $\lambda_i\in\mathbb{R}^{n_x}$ are the dual variables, $\ell_i:\mathbb{R}^{n_x}\times\mathbb{R}^{n_u}\times\mathbb{R}^{n_d}\rightarrow\mathbb{R}$ are the stage cost functions, $f_i:\mathbb{R}^{n_x}\times\mathbb{R}^{n_u}\times\mathbb{R}^{n_d}\rightarrow\mathbb{R}^{n_x}$ are the dynamic mapping functions, $\ell_N:\mathbb{R}^{n_x}\times\mathbb{R}^{n_d}\rightarrow\mathbb{R}$ is the final cost function, and the symbol $|$ is used to denote the associated dual variables. The initial state constraint is enforced with the initial state mapping $T\in\mathbb{R}^{n_0\times n_x}$ and parameter $d_{-1}\in\mathbb{R}^{n_0}$. We let $x_{-1},u_{-1},u_N,\lambda_N$ be empty vectors (for convenience), and define $z_{i}:=[x_i,u_i]$, $w_i:=[z_i;\lambda_i]$, and $\xi_i:=[w_i;d_i]$ for $i\in\mathbb{I}_{[-1,N]}$, and we use the syntax $v_{a:b}:=[v_{a};v_{a+1};\cdots;v_{b}]$ for $v=x,u,\lambda,d,z,w,\xi$. The DO problem \eqref{eqn:mpc} is a parametric nonlinear program that we denote as $P_{0:N}(d_{-1:N})$. We assume that all functions are twice continuously differentiable and potentially nonconvex. Typical MPC problems are formulated with $T=I$ and typical MHE problems are formulated with an empty matrix $T\in\mathbb{R}^{0\times n_x}$ (i.e., the initial constraint is not enforced). State-output mappings encountered in MHE problems are assumed to be embedded within the stage cost functions. 

In this paper we study a property of DO problems that is known as {\it exponential decay of sensitivity} (EDS) \citep{na2020exponential,shin2021exponential}. The property indicates that the sensitivity of the solution at stage $i$ against a data perturbation at stage $j$ decays exponentially with $|i-j|$. This property helps understand how different data perturbations (e.g., disturbances or changes in set-points, initial conditions, and terminal penalties) propagate along with the time horizon. Moreover, EDS has been shown to be essential in constructing efficient discretization schemes for continuous-time DO formulations \citep{shin2020diffusing,grune2020exponential} and in establishing convergence of algorithms \citep{na2020overlapping,na2020superconvergence}.

Building upon our previous results \citep{shin2021exponential} we show that, under uniform boundedness of the Lagrangian Hessian (uBLH), a uniform second order sufficiency condition (uSOSC), and a uniform linear independence constraint qualification (uLICQ), the primal-dual solution of the DO problem at a given stage decays exponentially with the distance to the stage at which a data perturbation is introduced.  In particular, given base data $d_{-1:N}^\star$ and associated primal-dual solution trajectory $w^\star_{-1:N}$ (at which uBLH, uSOSC, and uLICQ are satisfied), there exist uniform constants $\Upsilon>0$ and $\rho\in(0,1)$ and a neighborhood $\mathbb{D}^\star_{-1:N}$ of $d_{-1:N}^\star$ such that the following holds for any $d_{-1:N},d'_{-1:N}\in\mathbb{D}^\star_{-1:N}$: 
\begin{align}\label{eqn:eds}
  \| w^\dag_i(d_{-1:N})-w^\dag_i(d'_{-1:N})  \|  \leq  {\sum_{j=-1}^{N}}\Upsilon\rho^{|i-j|}\| d_j-d'_j\|,
\end{align}
where $w^\dag_i(\cdot)$ is the primal-dual solution mapping at stage $i$. That is, the sensitivity $\Upsilon\rho^{|i-j|}$ of the solution at stage $i$ against a data perturbation at stage $j$ decays exponentially with respect to the distance $|i-j|$. Here, it is important that $(\Upsilon,\rho)$ are uniform constants (independent of horizon length $N$); this allows us to maintain $(\Upsilon,\rho)$ unchanged even if the horizon length becomes indefinitely long (e.g., when approaching an infinite horizon). Our key finding is that uBLH, uLICQ, and uSOSC can be obtained under uniform controllability/observability and under uniformly bounded system matrices (standard assumption). This result thus allows to establish EDS directly from fundamental system-theoretic properties.

In summary, our main contribution is showing that controllability and observability provide sufficient conditions for EDS. This result sheds light on how system-theoretic properties influence the propagation of perturbations along the solution trajectory. EDS for continuous-time, linear-quadratic MPC has been established under stablizability and detectability in \citet{grune2019sensitivity,grune2020exponential,grune2020abstract}. EDS has been established for discrete-time, nonlinear MPC under uniform SOSC and controllability in \citep{na2020exponential}; the proof of this result uses a Riccati recursion representation of the optimality conditions of the DO problem. The proof that we present here is more compact and general; specifically, our approach is based on a graph-theoretic analysis of the optimality conditions. This general setting allows us to establish EDS for discrete-time, nonlinear MPC and MHE problems. 

{\it Basic Notation:} The set of real numbers and the set of integers are denoted by $\mathbb{R}$ and $\mathbb{I}$, respectively, and we define $\mathbb{I}_{A}:=\mathbb{I}\cap A$, where $A$ is a set; $\mathbb{I}_{>0}:=\mathbb{I}\cap(0,\infty)$; $\mathbb{I}_{\geq 0}:=\mathbb{I}\cap[0,\infty)$. We consider vectors always as column vectors. We use the syntax: $\{M_i\}_{i=1}^N:=[M_{1}; \cdots; M_{N}]:=[M_1^\top\,\cdots\,M_n^\top]^\top$. Furthermore, $v[i]$ denotes the $i$-th component of $v$. {For a function $\phi:\mathbb{R}^{n}\rightarrow \mathbb{R}$ and variable vectors $y\in\mathbb{R}^{p}$, $z\in\mathbb{R}^{q}$, $\nabla^2_{yz}\phi(x):=\{(\{\frac{\partial^2}{\partial y[i]\partial z[j]}\phi(x)\}_{j=1}^q)^\top\}_{i=1}^p$}. For a vector function $\varphi:\mathbb{R}^{n}\rightarrow \mathbb{R}^m$ and a variable vector $w\in\mathbb{R}^{s}$, $\nabla_{w}\varphi(x):=\{(\{\frac{\partial}{\partial w[j]}\varphi(x)[i]\}_{j=1}^s)^\top\}_{i=1}^m$. Vector 2-norms and induced 2-norms of matrices are denoted by $\|\cdot\|$. For matrices $A$ and $B$ with the same dimensions, $A\succ(\succeq) B$ indicates that $A-B$ is positive (semi) definite. We use the convention: if $m$ or $n$ is zero, $\mathbb{R}^{m\times n}$ is a singleton only containing the $m\times n$ null matrix.

\section{Main Results}\label{sec:main}
\subsection{Exponential Decay of Sensitivity}
In this section, we establish EDS \eqref{eqn:eds} for $P_{0:N}(\cdot)$. We first formally define sufficient conditions for EDS to hold (uBLH, uSOSC, and uLICQ). We begin by defining the notion of uniformly bounded quantities.
\begin{defn}[Uniform Bounds]
  A set $\{A_i\}_{i\in\mathcal{A}}$ is called $L$-uniformly bounded above (below) if there exists a uniform constant $L\in\mathbb{R}$ (independent of $N$) such that $\| A_i \| \leq (\geq) L$ holds for any $i\in\mathcal{A}$.
\end{defn}
We say that a set of a quantity is uniformly bounded above (below) if there exists a uniform constant $L$ such that the set is $L$-uniformly bounded above (below). Also, we will write that some quantity $a$ is uniformly bounded above (below) if $\{a\}$ is uniformly bounded above (below). 

The Lagrangian function of $P_{0:N}(d_{-1:N})$ is defined as
\begin{align*}
\mathcal{L}_{0:N}(w_{-1:N};d_{-1:N}):=\sum_{i=0}^N\mathcal{L}_{i}(z_i,\lambda_{i-1:i};d_i),
\end{align*}
where:
\begin{align*}
  \mathcal{L}_{i}(z_i,\lambda_{i-1:i};d_i)&:=\ell_i(z_i;d_i)-\lambda_{i-1}^\top x_i + \lambda_i^\top f_i(z_i;d_i)\\
  \mathcal{L}_N(x_N,\lambda_{N-1};d_N)&:=\ell_N(x_N;d_N) - \lambda^\top_{N-1} x_N.
\end{align*}

\begin{defn}[uBLH]
  Given $d_{-1:N}^\star$ and the solution $w_{-1:N}^\star$ of $P_{0:N}(d_{-1:N}^\star)$, $L$-uBLH holds if:
  \begin{align}\label{eqn:ublh}
    \Vert\nabla^2_{w_{-1:N}\xi_{-1:N}}\mathcal{L}_{0:N}(w_{-1:N}^\star;d_{-1:N}^\star)\Vert\leq L,
  \end{align}
  with uniform constant $L<\infty$.
\end{defn}

The primal Hessian $\bH_{0:N}$ of the Lagrangian and the constraint Jacobian $\bJ_{0:N}$ are:
\begin{subequations}\label{eqn:HJ}
\begin{align}
  \bH_{0:N}&:=\nabla^2_{z_{0:N},z_{0:N}}\mathcal{L}_{0:N}(w_{-1:N}^\star;d_{-1:N}^\star)\\
  \bJ_{0:N}&:=\nabla_{z_{0:N}}c_{-1:N-1}(z^\star_{0:N};d_{-1:N}^\star),
\end{align}
\end{subequations}
where $c_{-1:N-1}(\cdot)$ is the constraint function for $P_{0:N}(\cdot)$; {that is, $$c_{-1:N-1}(z_{0:N};d_{-1:N}):= \begin{bmatrix}Tx_{0}-d_{-1}\\x_{1}-f_1(z_0;d_0)\\\cdots\\x_{N}-f_{N-1}(z_{N-1};d_{N-1})\end{bmatrix}.$$}

\begin{defn}[uSOSC]
  Given $d_{-1:N}^\star$ and the solution $w_{-1:N}^\star$ of $P_{0:N}(d_{-1:N}^\star)$, $\gamma$-uSOSC holds if:
  \begin{align}\label{eqn:usosc}
    ReH(\bH_{0:N},\bJ_{0:N})\succeq \gamma I,
  \end{align}
  with uniform constant $\gamma>0$.
\end{defn}
Here, $ReH(\bH_{0:N},\bJ_{0:N}):=Z^\top \bH_{0:N}Z$ is the reduced Hessian and $Z$ is a null-space matrix of $\bJ_{0:N}$. 
\begin{defn}[uLICQ]
  Given $d_{-1:N}^\star$ and the primal-dual solution $w_{-1:N}^\star$ of $P_{0:N}(d_{-1:N}^\star)$, $\beta$-uLICQ holds if:
  \begin{align}\label{eqn:ulicq}
  \bJ_{0:N}\bJ_{0:N}^\top \succeq \beta I
  \end{align}
  with uniform constant $\beta>0$.
\end{defn}
Note that uSOSC assumes that the smallest eigenvalue of the reduced Hessian is uniformly bounded below by $\gamma$, while uLICQ assumes that the smallest non-trivial singular value of the Jacobian is uniformly bounded below by ${\beta}^{1/2}$. Thus, these are strengthened versions of SOSC and LICQ.
We require uSOSC and uLICQ because, under SOSC and LICQ, the smallest eigenvalue of reduced Hessian or the smallest non-trivial singular value of the Jacobian may become arbitrarily close to $0$ as the horizon length $N$ is extended (e.g., see \citet[Example 4.18]{shin2021exponential}). Under uSOSC and uLICQ, on the other hand, the lower bounds are independent of $N$.

\begin{assum}\label{ass:regularity}
  Given twice continuously differentiable   \linebreak functions $\{\ell_i(\cdot)\}_{i=0}^{N}$, $\{f_i(\cdot)\}_{i=0}^{N-1}$ and base data $d^\star_{-1:N}$, there exists a primal-dual solution $w_{-1:N}^\star$ of $P_{0:N}(d_{-1:N}^\star)$ at which $L$-uBLH, $\gamma$-uSOSC, and $\beta$-uLICQ are satisfied.
\end{assum}
The following lemma is a well-known characterization of  solution mappings of parametric nonlinear programs (NLPs) \citep{robinson1980strongly,dontchev2009implicit}.
\begin{lem}\label{lem:fun}
  Under Assumption \ref{ass:regularity}, there exist neighborhoods $\mathbb{D}^\star_{-1:N}$ of $d^\star_{-1:N}$ and $\mathbb{W}^\star_{-1:N}$ of $w^\star_{-1:N}$ and continuous $w^\dag_{-1:N}:\mathbb{D}^\star_{-1:N}\rightarrow \mathbb{W}^\star_{-1:N}$ such that for any $d_{-1:N}\in\mathbb{D}^\star_{-1:N}$, $w^\dag_{-1:N}(d_{-1:N})$ is a local solution of $P_{0:N}(d_{-1:N})$.
\end{lem}
\begin{pf}
  From \citet[Lemma 3.3]{shin2021exponential}.
\qed\end{pf}
We can thus see that there exists a well-defined solution mapping $w^\dag_{-1:N}(\cdot)$ around the neighborhood of $d_{-1:N}^\star$. We now study {\it stage-wise} solution sensitivity by characterizing the dependence of $w^\dag_i(\cdot)$ on the data $d_{-1:N}$. 

\begin{thm}\label{thm:eds}
  Under Assumption \ref{ass:regularity}, there exist uniform constants $\Upsilon>0$ and $\rho\in(0,1)$ (functions of $L,\gamma,\beta$) and neighborhoods $\mathbb{D}^\star_{-1:N}$ of $d_{-1:N}^\star$ and $\mathbb{W}^\star_{-1:N}$ of $w_{-1:N}^\star$ such that \eqref{eqn:eds} holds for any $d_{-1:N},d'_{-1:N}\in\mathbb{D}^\star_{-1:N}$ and $i\in \mathbb{I}_{[-1,N]}$.
\end{thm}
\begin{pf}
  We observe that $P_{0:N}(\cdot)$ is graph-structured (induced by $\mathcal{G}_N=(\mathcal{V}_N,\mathcal{E}_N)$, where $\mathcal{V}_N=\{-1,0,\cdots,N\}$ and $\mathcal{E}_N=\{\{-1,0\},\{0,1\},\cdots,\{N-1,N\}$), and the maximum graph degree $D=2$. From uBLH, uLICQ, and uSOSC, one can see that assumptions in \citet[Theorem 4.9]{shin2021exponential} are satisfied. This implies that the singular values of $\nabla^2_{w_{0:N}w_{0:N}}\mathcal{L}_{0:N}(w^\star_{-1:N};d^\star_{-1:N})$ are uniformly upper and lower bounded and those of $\nabla^2_{w_{0:N}d_{-1:N}}\mathcal{L}_{0:N}(w^\star_{-1:N};d^\star_{-1:N})$ are uniformly upper bounded (uniform constants given by functions of $L,\beta,\gamma$; see \citet[Equation (4.15)]{shin2021exponential}). We then apply \citet[Theorem 3.5]{shin2021exponential} to obtain $\Upsilon>0$ and $\rho\in(0,1)$ as functions of the upper and lower bounds of the singular values (see \citet[Equation (3.17)]{shin2021exponential}). This allows expressing $\Upsilon,\rho$ as functions of $L, \beta, \gamma$.
  \qed\end{pf}

Theorem \ref{thm:eds} establishes EDS under the regularity conditions of Assumption \ref{ass:regularity}. It is important that $\Upsilon,\rho$ can be determined solely in terms of $L,\gamma,\beta$ (and do not depend on the horizon length $N$). Practical DO problems typically have additional equality/inequality constraints that are not considered in \eqref{eqn:mpc}. Thus, Theorem \ref{thm:eds} may not be directly applicable to those problems. However, the results in \citet{shin2021exponential} are applicable to such problems as long as the DO problem is a graph-structured NLP. Specifically, under uniformly strong SOSC and uLICQ, we can establish EDS using \citet[Theorem 3.5, 4.9]{shin2021exponential}. {The graph structure breaks when there exist globally coupled variables; typical MPC and MHE problems do not have such variables, but parameter estimation problems may have such variables. Specifically, in the presence of globally coupled variables, the graph distance between any pair of stages is not greater than two.}

\subsection{Regularity from System-Theoretic Properties}
Although uSOSC and uLICQ are standard notions of NLP solution regularity, they are not intuitive notions from a system-theoretic perspective. However, we now show that uSOSC and uLICQ can be obtained from uniform controllability and observability.  We begin by defining:
\begin{align*}
  Q_i&:=\nabla^2_{x_ix_i} \mathcal{L}_i(z^\star_i,\lambda^\star_{i-1:i};d^\star_i)
  &&R_i:=\nabla^2_{u_iu_i} \mathcal{L}_i(z^\star_i,\lambda^\star_{i-1:i};d^\star_i)\\
  S_i&:=\nabla^2_{x_iu_i} \mathcal{L}_i(z^\star_i,\lambda^\star_{i-1:i};d^\star_i)
  &&E_i:=\nabla^2_{x_id_i} \mathcal{L}_i(z^\star_i,\lambda^\star_{i-1:i};d^\star_i)\\
  F_i&:=\nabla^2_{u_id_i} \mathcal{L}_i(z^\star_i,\lambda^\star_{i-1:i};d^\star_i)
  &&A_i:=\nabla_{x_i}f_i(z^\star_i;d^\star_i)\\
  B_i&:=\nabla_{u_i}f_i(z^\star_i;d^\star_i)
  &&G_i:=\nabla_{d_i}f_i(z^\star_i;d^\star_i).
\end{align*}

\begin{defn}
  $(\{A_i\}_{i=1}^{N-1},\{B_i\}_{i=0}^{N-1})$ is $(N_c,\beta_c)$-uniformly controllable with $N_c\in\mathbb{I}_{\geq 0}$ and $\beta_c>0$ (independent of $N$) if, for any $i,j\in\mathbb{I}_{[0,N-1]}$ with $|i-j|\geq N_c$, $\mathcal{C}_{i:j} \mathcal{C}_{i:j}^\top\succeq \beta_c I$ holds, where
  \begin{align*}
    \mathcal{C}_{i:j}:=
    \begin{bmatrix}
      A_{i+1:j} B_{i} & \cdots & A_{j} B_{j-1} & B_{j}
    \end{bmatrix}.
  \end{align*}
\end{defn}
\begin{defn}
  $(\{A_i\}_{i=0}^{N-1},\{Q_i\}_{i=0}^{N})$ is $(N_o,\gamma_o)$-uniformly observable with $N_o\in\mathbb{I}_{\geq 0}$ and $\gamma_o>0$ (independent of $N$) if for any $i,j\in\mathbb{I}_{[0,N-1]}$ with $|i-j|\geq N_o$, $\mathcal{O}^\top_{i:j} \mathcal{O}_{i:j}\succeq \gamma_o I$ holds, where
  \begin{align*}
    \mathcal{O}_{i:j}:=
    \begin{bmatrix}
      Q_{j}A_{i:j-1} \\ \ddots \\  Q_{i+1}A_{i} \\ Q_{i}
    \end{bmatrix}.
  \end{align*}
\end{defn}
Here $$A_{a:b}:=
\begin{cases}
  A_{b}A_{b-1}\cdots A_{a+1}A_{a},\text{ if }a\leq b\\
  A_{b}A_{b+1}\cdots A_{a-1}A_{a},\text{ otherwise.}
\end{cases}
$$
Note that uniform controllability and observability are stronger versions of their standard counterparts. One can establish the following duality between uniform controllability and observability.
\begin{prop}\label{prop:dual}
  $(\{A_i\}_{i=1}^{N},\{B_i\}_{i=0}^{N})$ is $(N_0,\alpha_0)$-uniformly controllable if and only if $(\{A^\top_i\}_{i=N}^1,\{B^\top_i\}_{i=N}^0)$ is $(N_0,\alpha_0)$-uniformly observable (here, note that the orders of sequences $\{A^\top_i\}_{i=N}^1,\{B^\top_i\}_{i=N}^0$ are inverted).  
\end{prop}
\begin{pf}
  {The proof is straightforward and thus omitted.} \qed
\end{pf}


The following technical lemma is needed to show that uniform controllability implies uLICQ.
\begin{lem}\label{lem:blk}
  Consider a block row/column operator $U$ with $L$-uniformly bounded above block $V$ of the form:
  \begin{align*}
    U:=
    \begin{bmatrix}
      I\\
      V& I\\
      &&\ddots \\
      &&&I
    \end{bmatrix},
          \begin{bmatrix}
            I&&&V\\
            & I\\
            &&\ddots \\
            &&&I
          \end{bmatrix}.
  \end{align*}
  We have that $U,U^{-1}$ are $(L+1)$-uniformly bounded above.
\end{lem}
\begin{pf}
  {The proof is straightforward and thus omitted.} \qed
\end{pf}
  We now show that uniform controllability implies uLICQ.
\begin{lem}\label{lem:ulicq}
  $K$-uniform upper boundedness of $\{A_i\}_{i=0}^{N-1}$ and $\{B_i\}_{i=0}^{N-1}$, $TT^\top\succeq \delta I$ for uniformly lower bounded  \linebreak $\delta>0$, and $(N_c,\beta_c)$-uniform controllability of   \linebreak $(\{A_i\}_{i=1}^{N-1},\{B_i\}_{i=0}^{N-1})$ implies \eqref{eqn:ulicq}, where $\beta>0$ is a function of $K,\delta,N_c,\beta_c$.
\end{lem}
\begin{pf}
  The Jacobian $\bJ_{0:N}$ has the following form:
  \begin{align*}\small
    \bJ_{0:N}=
    \begin{bmatrix}
      T\\
      -A_0&-B_{0}&I\\
      &&\ddots\\
      &&-A_{N-2}&-B_{N-2}&I\\
      &&&&-A_{N-1}&-B_{N-1}& I
    \end{bmatrix}.
  \end{align*}
  By inspecting the block structure of $\bJ_{0:N}$ and \citet[Lemma 4.15]{shin2021exponential}, one can see that it suffices to show that the smallest non-trivial singular value of 
  \begin{align}\label{eqn:mat-1}
    \begin{bmatrix}
      S\\
      -A_i&-B_{i}&I\\
      &&\ddots\\
      &&-A_{j-1}&-B_{j-1}&I\\
      &&&&-A_{j}&-B_{j}
    \end{bmatrix}
  \end{align}
  is ${\beta}^{1/2}$-uniformly bounded below for $S=T$ or $I$ and for any $i,j\in\mathbb{I}_{[0,N-1]}$ with $ N_c\leq |i-j|\leq  2N_c$, where $0<\beta\leq 1$ is a function of $K,\delta,N_c,\beta_c$. This follows from the observation that one can always partition $\mathbb{I}_{[0,N-1]}$ into a family of blocks with size between $N_c$ and $2N_c$. For now, we assume $S=I$. By applying a set of suitable block row and column operations (in particular, first apply block row operations to eliminate $A_i,\cdots,A_j$, and then apply block column operations to eliminate $-B_i,\cdots,-A_{i:j-1}B_{j-2}$) and permutations, one can obtain the following:
  \begin{align}\label{eqn:mat-3}
    \begin{bmatrix}
      I\\
      &-A_{i+1:j}B_i&\cdots&-A_{j}B_{j-1}&-B_j
    \end{bmatrix}.
  \end{align}
  The lower-right blocks constitute the controllability matrix $\mathcal{C}_{i:j}$; from uniform controllability, the smallest non-trivial singular value of the matrix in \eqref{eqn:mat-3} is uniformly lower bounded by $\min(1,\beta_c^{1/2})$. Here, we have applied block-row and block-column operations as the ones that appear in Lemma \ref{lem:blk} (each multiplied block is uniformly bounded above due to $K$-uniform boundedness of $\{A_i\}_{i=0}^{N-1}$ and $\{B_i\}_{i=0}^{N-1}$). Also, we have applied such operations only uniformly bounded many times (the number of operations is independent of $N$ since the number of blocks in the matrix in \eqref{eqn:mat-1} is bounded by $4(2N_c+1)(N_c+1)$, which is uniformly bounded above). We thus have that the smallest non-trivial singular value of the matrix in \eqref{eqn:mat-1} is uniformly lower bounded with uniform constant ${\beta_0}^{1/2}$, and $\beta_0>0$ is given by a function of $K,N_c,\beta_c$. Now we consider the $S=T$ case. One can observe that the smallest non-trivial singular value of the matrix in \eqref{eqn:mat-1} with $S=T$ is lower bounded by that with $S=[\widetilde{T};T]$ (here, $\widetilde{T}^\top$ is a null space matrix of $T$); and again, it is lower bounded by $\delta^{1/2}$ times that with $S=I$. We thus have that the smallest non-trivial singular value of the matrix in \eqref{eqn:mat-1} with $S=T$ is uniformly lower bounded by $\beta_0^{1/2}\delta^{1/2}$. Therefore, the smallest non-trivial singular values of the matrices in \eqref{eqn:mat-1} with $S=I$ or $T$ are $\beta^{1/2}$-uniformly lower bounded for any $i,j\in\mathbb{I}_{[0,N-1]}$ with $ N_c\leq |i-j|\leq  2N_c$, where $\beta:=\min(\beta_0,\delta\beta_0,1)$. Thus, by \citet[Lemma 4.15]{shin2021exponential} we have \eqref{eqn:ulicq}. 
  \qed\end{pf}

If $T\in\mathbb{R}^{0\times n_x}$, the assumption $TT^\top\succeq \delta I$ for uniformly lower bounded $\delta>0$ holds for an arbitrary $\delta>0$ due to the convention introduced in the {\it Notation} in Section \ref{sec:intro}.

We now show that uniform observability implies uSOSC.
\begin{lem}\label{lem:usosc}
  $K$-uniform boundedness of $\{A_i\}_{i=0}^{N-1}$, $\{B_i\}_{i=0}^{N-1}$, and $\{Q_i\}_{i=0}^{N}$, $Q_i\succeq 0$, $S_i=0$, $R_i\succeq rI$ ($r>0$ is independent of $N$), and $(N_o,\gamma_o)$-uniform observability of $(\{A_i\}_{i=0}^{N-1},\{Q_i\}_{i=0}^{N})$ implies \eqref{eqn:usosc}, where $\gamma>0$ is a function of $K,N_o,\gamma_o,r$.
\end{lem}
\begin{pf}
  The primal Hessian $\bH_{0:N}$ has the following form:
  \begin{align*}
    \bH_{0:N}
    \begin{bmatrix}
    Q_0\\
    &R_0\\
    &&\ddots\\
    &&&Q_{N-1}\\
    &&&&R_{N-1}\\
    &&&&&Q_N
  \end{bmatrix}.
  \end{align*}
  By inspecting the block structure of $\bH_{0:N}$ and $\bJ_{0:N}$ and  \citet[Lemma 4.14]{shin2021exponential}, one can observe that it suffices to show that: first, 
  \begin{align}\label{eqn:mat-5}
ReH\left(\begin{bmatrix}\bQ_{i:j}\\&\bR_{i:j-1}\end{bmatrix},\begin{bmatrix}\bA_{i:j}& \bB_{i:j-1}\end{bmatrix}\right)
  \end{align}
  has $\gamma$-uniformly lower bounded smallest eigenvalue with $\gamma>0$ for any $i,j\in\mathbb{I}_{[0,N-1]}$ with $N_o\leq |i-j| \leq 2N_o$, where:
  \begin{small}
    \begin{align*}
      & \bA_{i:j}:=
        \begin{bmatrix}
          -A_i & I              \\
          & \ddots & \ddots         \\
          &  &   -A_{j-1} & I
        \end{bmatrix},
                            \bB_{i:j-1}:=
                            \begin{bmatrix}
                              -B_i          \\
                              &    & \ddots    \\
                              &    &    & -B_{j-1}
                            \end{bmatrix} \\
      & \bQ_{i:j}:=
        \begin{bmatrix}
          Q_i                   \\
          & Q_{i+1}              \\
          &    & \ddots    \\
          &    &    & Q_{j}
        \end{bmatrix},
                      \bR_{i:j-1}:=
                      \begin{bmatrix}
                        R_i        \\
                        & R_{i+1}              \\
                        &    & \ddots    \\
                        &    &    & R_{j-1}
                      \end{bmatrix},
    \end{align*}
  \end{small}
    and second, $R_i\succeq \gamma I$ for any $i\in\mathbb{I}_{[0,N-1]}$. This follows from the observation that one can always partition $\mathbb{I}_{[0,N-1]}$ into a family of blocks with size between $N_c$ and $2N_c$. We consider $\bx_{i:j},\bu_{i:j-1}$ such that $\bA_{i:j}\bx_{i:j}+\bB_{i:j-1}\bu_{i:j-1} = 0$ holds. By uniform positive definiteness of $\{R_i\}_{i=0}^{N-1}$ and uniform boundedness of $\{B_{i}\}_{i=0}^{N-1}$, for $\kappa:=r/2K^2$, we have
  \begin{align*}
    \frac{1}{2}\bu_{i:j-1}^\top\bR_{i:j-1}\bu_{i:j-1} 
                                                      &\geq \kappa \bu_{i:j-1}^\top \bB_{i:j-1}^\top \bB_{i:j-1} \bu_{i:j-1}\\
    &= \kappa \bx_{i:j}^\top \bA_{i:j}^\top \bA_{i:j} \bx_{i:j},
  \end{align*}
  where the equality follows from $\bA_{i:j}\bx_{i:j}+\bB_{i:j-1}\bu_{i:j-1} = 0$. Furthermore, from $\bQ_{i:j}\succeq 0$, we have that:
  \begin{align*}
    \bx_{i:j}^\top\bQ_{i:j}^2\bx_{i:j} = (\bQ_{i:j}^{1/2}\bx_{i:j})^\top\bQ_{i:j}(\bQ_{i:j}^{1/2}\bx_{i:j})\leq K \bx_{i:j}^\top \bQ_{i:j} \bx_{i:j},
  \end{align*}
  where the inequality follows from that the largest eigenvalue of $\bQ_{i:j}$ is bounded by $K$.
  Thus, $\bx_{i:j}^\top \bQ_{i:j} \bx_{i:j}+\bu_{i:j-1}^\top\bR_{i:j-1}\bu_{i:j-1}$ is not less than:
  \begin{align*}
     \min(1/K,\kappa)  \bx_{i:j}^\top
      \begin{bmatrix}
        \bQ_{i:j}& \bA_{i:j}^\top
    \end{bmatrix}
    \begin{bmatrix}
      \bQ_{i:j}\\ \bA_{i:j}
    \end{bmatrix} \bx_{i:j} + \frac{r}{2} \|\bu_{i:j-1}\|^2.
  \end{align*}
  Observe that $\begin{bmatrix} \bQ_{i:j}& \bA_{i:j}^\top \end{bmatrix}$ can be permuted to:
  \begin{align}\label{eqn:mat-6}
    \begin{bmatrix}
      Q_{j}&I\\
      &-A_{j-1}&Q_{j-1}&I\\
      &&&\ddots\\
      &&&-A^\top_{i+1}&Q_{i+1}& I\\
      &&&&&-A^\top_{i}&Q_i
    \end{bmatrix}.
  \end{align}
We apply block row and column operations (as those of Lemma \ref{lem:fun}) uniformly bounded many times to obtain:
  \begin{align}\label{eqn:mat-7}
    \begin{bmatrix}
      I\\
      &A^\top_{j-1:i}Q_{j}&\cdots&A^\top_{i}Q_{i+1}& Q_i\\
    \end{bmatrix}.
  \end{align}
  From Proposition \ref{prop:dual} and $(N_o,\gamma_o)$-uniform observability of $(\{A_i\}_{i=0}^{N-1},\{Q_i\}_{i=0}^{N})$, we have that
\begin{equation}  
  (\{A^\top_i\}_{i=N-1}^{0},\{Q_i\}_{i=N}^{0})
\end{equation}  
   is $(N_o,\gamma_o)$-uniformly controllable. We thus have that the matrix in \eqref{eqn:mat-7} has $\min(1,\gamma_o)$-uniformly lower bounded smallest non-trivial singular value. This implies that the smallest non-trivial singular value of the matrix in \eqref{eqn:mat-6} is uniformly lower bounded by $\gamma'$, where $\gamma'$ is given by a function of $K$, $N_o$, $\gamma_o$. Therefore, we have that:
   $\bx_{i:j}^\top \bQ_{i:j} \bx_{i:j}+\bu_{i:j-1}^\top\bR_{i:j-1}\bu_{i:j-1}
    \geq \gamma (\| \bx_{i:j}\|^2+ \| \bu_{i:j-1}\|^2)$,
  where $\gamma:=\min(\gamma'/K,\kappa\gamma',r/2)$. One can observe that $R_i\succeq \gamma I$ for any $i\in\mathbb{I}_{[0,N-1]}$. Consequently, the smallest eigenvalues of the matrix in \eqref{eqn:mat-5} for any $i,j\in\mathbb{I}_{[0,N-1]}$ with $ N_o\leq |i-j|\leq  2N_o$ and $R_i$ are $\gamma$-uniformly lower bounded. Thus, by \citet[Lemma 4.14]{shin2021exponential}, \eqref{eqn:usosc} holds. One can confirm that $\gamma$ is a function of $K,N_o,\gamma_o,r$.
  \qed\end{pf}

Finally, we show that uniform boundedness of system matrices implies uBLH.  

\begin{lem}\label{lem:ub}
  If $\{Q_i\}_{i=0}^{N}$, $\{R_i\}_{i=0}^{N-1}$, $\{S_i\}_{i=0}^{N-1}$, $\{A_i\}_{i=0}^{N-1}$, $\{B_i\}_{i=0}^{N-1}$, $\{E_i\}_{i=0}^{N-1}$, $\{F_i\}_{i=0}^{N-1}$, $\{G_i\}_{i=0}^{N-1}$, and $T$ are $K$-uniformly bounded above, \eqref{eqn:ublh} holds, where $L<\infty$ is a function of $K$.
\end{lem}
\begin{pf}
  Uniform boundedness of the system matrices implies that for any $i,j\in\mathbb{I}_{[-1,N]}$, $\nabla^2_{w_i\xi_i}\mathcal{L}_{0:N}(w^\star_{-1:N};d^\star_{-1:N})$ are uniformly bounded above by $4K$ (\citet[Lemma 4.6]{shin2021exponential}). Furthermore, by inspecting the problem structure, we can see that $\Vert\nabla^2_{w_i\xi_j}\mathcal{L}_{0:N}(w^\star_{-1:N};d^\star_{-1:N})\Vert\leq 1$ for $i\neq j$ (there is only one identity block). Thus, $\Vert\nabla^2_{w_i\xi_j}\mathcal{L}_{0:N}(w^\star_{-1:N};d^\star_{-1:N})\Vert \leq \max(4K,1)$. 
  By noting that the maximum graph degree $D=2$ and applying \citet[Lemma 4.5]{shin2021exponential}, we have that $\nabla^2_{w_{-1:N}\xi_{-1:N}}\mathcal{L}_{0:N}(w^\star_{-1:N};d^\star_{-1:N})$ is $4\max(4K,1)$-uniformly bounded above. We set $L:=4\max(4K,1)$.\qed
\end{pf}
We now state EDS in terms of uniformly bounded system matrices and uniform controllability/observability. 
\begin{assum}\label{ass:regularity-2}
  Given twice continuously differentiable functions $\{\ell_i(\cdot)\}_{i=0}^{N}$, $\{f_i(\cdot)\}_{i=0}^{N-1}$, and data $d^\star_{-1:N}$, there exists a primal-dual solution $w_{-1:N}^\star$ of $P_{0:N}(d_{-1:N}^\star)$ at which the assumptions in Lemma \ref{lem:ulicq}, \ref{lem:usosc}, \ref{lem:ub} hold.
\end{assum}
\begin{cor}\label{cor:sys}
  Under Assumption \ref{ass:regularity-2}, there exist uniform constants $\Upsilon>0$ and $\rho\in(0,1)$ (functions of $K$, $r$, $N_c$, $\beta_c$, $N_o$, $\gamma_o$, $\delta$) and neighborhoods $\mathbb{D}^\star_{-1:N}$ of $d_{-1:N}^\star$ and $\mathbb{W}^\star_{-1:N}$ of $w_{-1:N}^\star$ such that \eqref{eqn:eds} holds for any $d_{-1:N},d'_{-1:N}\in\mathbb{D}^\star_{-1:N}$ and $i\in \mathbb{I}_{[-1,N]}$.
\end{cor}
\begin{pf}
 From Theorem \ref{thm:eds} and Lemma \ref{lem:ulicq}, \ref{lem:usosc}, \ref{lem:ub}.
\qed\end{pf}
\subsection{Time-Invariant Setting}\label{sec:ti}
Assume now that the system is time-invariant and focus on a region around a steady-state. A corollary of Theorem \ref{thm:eds} for such a setting is derived. We present this result since this setting has been of particular interest in the MPC literature. Consider a time-invariant system with a stage-cost function $\ell(\cdot)$, initial regularization function $\ell_b(\cdot)$, terminal cost function $\ell_f(\cdot)$, and dynamic mapping $f(\cdot)$. The DO problem is given by \eqref{eqn:mpc} with  $f_i(\cdot)=f(\cdot)$ for $i\in\mathbb{I}_{[0,N-1]}$, $\ell_i(\cdot)=\ell(\cdot)$ for $i\in\mathbb{I}_{[1,N-1]}$, $\ell_0(x,u)=\ell(x,u;d)+\ell_b(x;d)$, and $\ell_N(x;d)=\ell_f(x;d)$. The steady-state optimization problem is:
\begin{align}\label{eqn:ss}
  \min_{x,u}\;& \ell(x,u;d)\;  \st\;x=f(x,u;d) \quad|\quad\lambda.
\end{align}
For given $d^s$ and an associated primal-dual solution $w^s:=[x^s;u^s;\lambda^s]$ of \eqref{eqn:ss}, we define:
\begin{align*}
  Q&:=\nabla^2_{xx}\mathcal{L}(w^s;d^s),\, S:=\nabla^2_{xu}\mathcal{L}(w^s;d^s),\\
  R&:=\nabla^2_{uu}\mathcal{L}(w^s;d^s),\,
  A:=\nabla_{x}f(z^s;d^s),\, B:=\nabla_{u}f(z^s;d^s),
\end{align*}
where $\mathcal{L}(w;d):=f(z;d)-\lambda^\top x + \lambda^\top f(z;d)$; for the initial and terminal cost functions $\ell_b(\cdot)$ and $\ell_f(\cdot)$, we define:
\begin{align*}
  \lambda_b&:=\nabla_{x}\ell_b(x^s;d^s),\;Q_b:=\nabla^2_{xx}\ell_b(x^s;d^s)\\
  \lambda_f&:=\nabla_{x}\ell_f(x^s;d^s),\;Q_f:=\nabla^2_{xx}\ell_f(x^s;d^s).
\end{align*}
The quantities defined above ($Q$, $R$, etc.) are independent of $N$ since $w^s$ can be determined independently of $N$.

\begin{assum}\label{ass:regularity-3}
  Given twice continuously differentiable $\ell(\cdot)$, $\ell_b(\cdot)$, $\ell_f(\cdot)$, $f(\cdot)$, and data $d^s$, there exists a steady-state solution $w^s$, at which $Q_f\succeq Q\succeq 0$, $Q_b\succeq 0$, $S=0$, $R\succ 0$, $(A,B)$ controllable, $(A,Q)$ observable, $TT^\top \succ 0$, $\lambda_b+\lambda^s\in\text{Range}(T^\top)$ and $\lambda_f=\lambda^s$ hold.
\end{assum}
\begin{cor}\label{cor:ti}
  Under the time invariance setting and Assumption \ref{ass:regularity-3}, there exist uniform constants $\Upsilon>0$ and $\rho\in(0,1)$ such that the following holds: for any $N\in\mathbb{I}_{\geq 0}$, there exist neighborhoods $\mathbb{D}^s_{-1:N}$ of $d^s_{-1:N}:=[Tx^s;d^s;\cdots;d^s]$ and $\mathbb{W}^s_{-1:N}$ of $w_{-1:N}^s:=[\lambda^s_{-1};w^s;\cdots;w^s;x^s]$ such that \eqref{eqn:eds} holds for any $d_{-1:N},d'_{-1:N}\in\mathbb{D}^s_{-1:N}$, where $\lambda^s_{-1}$ is the solution of $T^\top\lambda^{s}_{-1} = \lambda_b + \lambda^s$.
\end{cor}
\begin{pf}
  From the existence (follows from $\ell_b+\ell^s\in\text{Range}(T^\top)$) and uniqueness (follows from $TT^\top \succ 0$) of the solution of $T^\top\lambda^{s}_{-1} = \lambda_b + \lambda^s$, we have well-defined $\lambda^s_{-1}$. From $T^\top\lambda^{s}_{-1} = \lambda_b + \lambda^s$ and the optimality of $w^s$ for \eqref{eqn:ss}, $w^s_{-1:N}$ satisfies the first-order optimality conditions for $P_{0:N}(d_{-1:N}^s)$. Furthermore, all the assumptions in Lemma \ref{lem:ub} are satisfied with some uniform constant $K$ because $\ell(\cdot)$, $\ell_b(\cdot)$, $\ell_f(\cdot)$, $f(\cdot)$, $T$, $w^s$, and $ d^s$ are independent of $N$; thus, by Lemma \ref{lem:ub}, we have \eqref{eqn:ublh} for a uniform constant $L<\infty$. Moreover, $TT^\top\succeq \delta I$ holds for some uniform constant $\delta>0$, and $R_i\succeq rI$ for $i\in\mathbb{I}_{[0,N-1]}$ with some uniform constant $r>0$, since $\ell(\cdot)$, $w^s$, $d^s,T$ are independent of $N$. Similarly, $(A,B)$ controllability implies $(N_c,\beta_c)$-uniform controllability of $(\{A_i\}_{i=1}^{N-1},\{B_i\}_{i=0}^{N-1})$ with some uniform constant $N_c,\beta_c$, and $(A,Q)$ observability implies $(N_o,\gamma_o)$-uniform observability of $(\{A_i\}_{i=0}^{N-1},\{Q_i\}_{i=0}^{N})$ for some uniform constants $N_o,\gamma_o$ (for now, we assume that $Q_b=0$ and $Q_f=Q$). From Lemma \ref{lem:ulicq}, \ref{lem:usosc}, we have \eqref{eqn:usosc} and \eqref{eqn:ulicq} for uniform $\beta,\gamma>0$. Now, observe that \eqref{eqn:usosc} for $Q_b=0$ and $Q_f=Q$ implies \eqref{eqn:usosc} for any $Q_b\succeq 0$ and $Q_f\succeq Q$; thus, we have \eqref{eqn:usosc} with uniform $\gamma>0$ for any $Q_b,Q_f$. Since the first and second order conditions of optimality and constraint qualifications are satisfied, $w^s_{-1:N}$ is a strict minimizer for $P_{0:N}(d_{-1:N})$. Since we have \eqref{eqn:ublh}, \eqref{eqn:usosc}, and \eqref{eqn:ulicq} with uniform $L,\gamma,\beta$, we have uBLH, uLICQ, and uSOSC at $(w^s_{-1:N},d^s_{-1:N})$. By applying Theorem \ref{thm:eds}, we can obtain \eqref{eqn:eds}. Lastly, since the parameters $K,r,N_c,\beta_c,N_o,\gamma_o$ are independent of $N$, so do $\Upsilon$ and $\rho$.
\end{pf}

Initial and terminal cost functions that satisfy Assumption \ref{ass:regularity-3} can be constructed as:
\begin{align*}
  \ell_b(x;d) &:=-((I-T^+T)\lambda^s)^\top x\\
  \ell_f(x;d) &:= ( x-x^s)^\top Q( x-x^s) + (\lambda^s )^\top x,
\end{align*}
where $(\cdot)^+$ is the pseudoinverse of the argument. One can observe that $\ell_b(\cdot)$ can be set to constantly zero if $T=I$.

\section{Numerical Results}\label{sec:num}
\begin{figure*}[t]
  \centering
  \begin{tabular}{cc}
    \includegraphics[width=.43\textwidth]{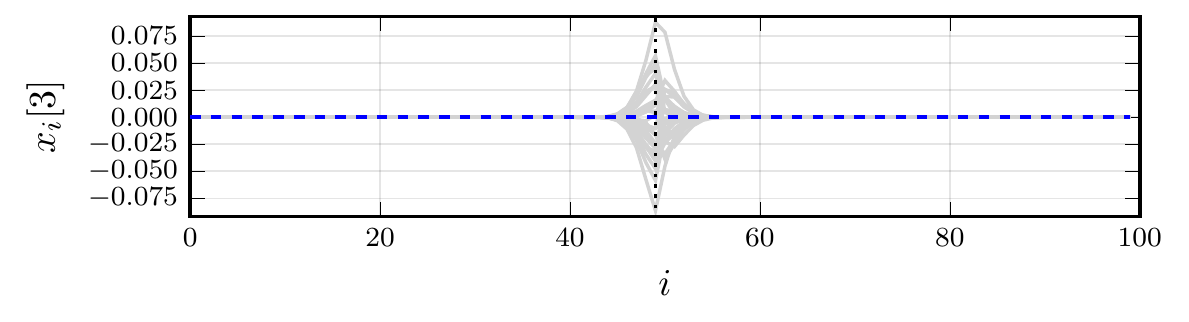}&
    \includegraphics[width=.43\textwidth]{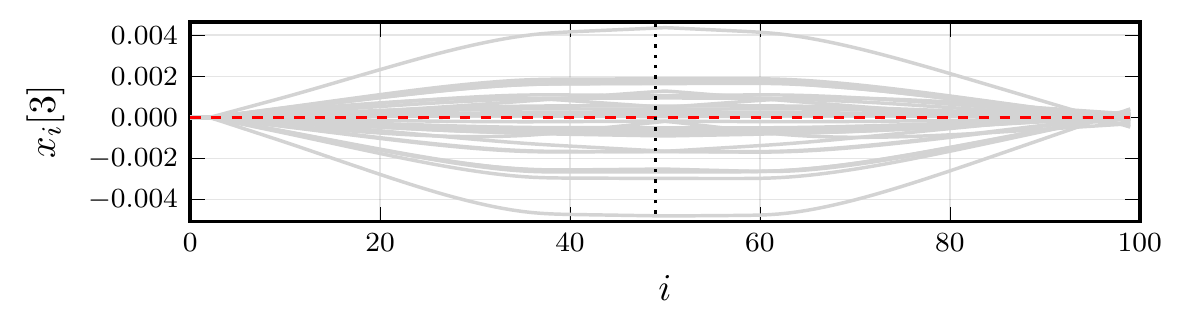}\\
    \includegraphics[width=.43\textwidth]{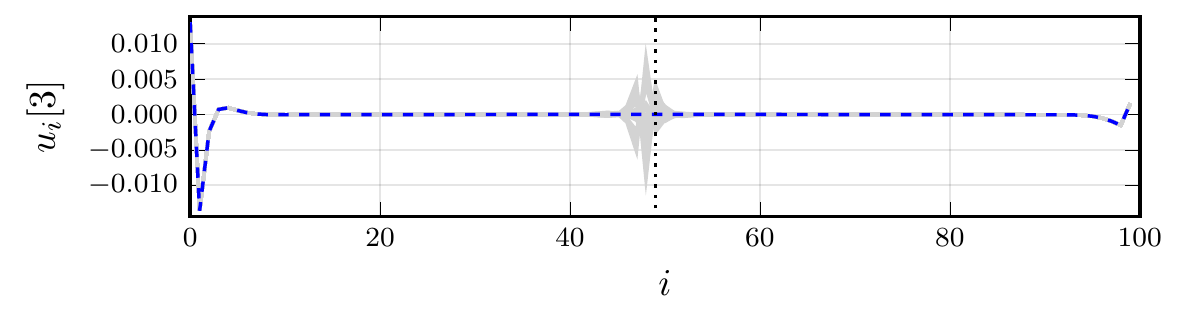}&
    \includegraphics[width=.43\textwidth]{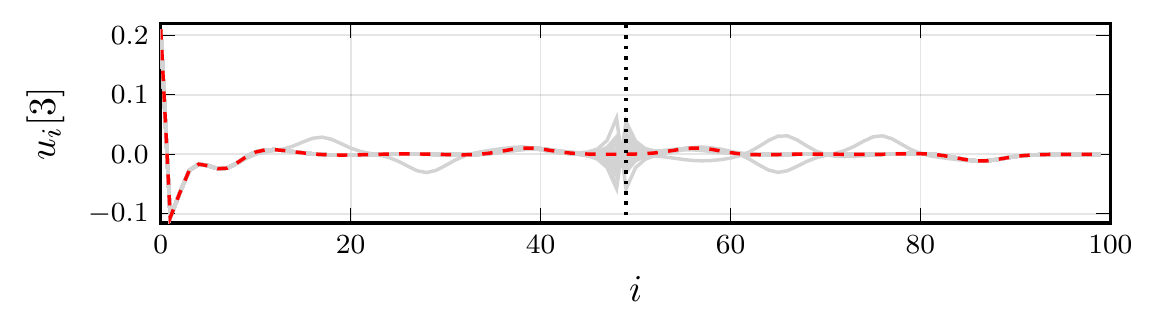}\\
    \includegraphics[width=.43\textwidth]{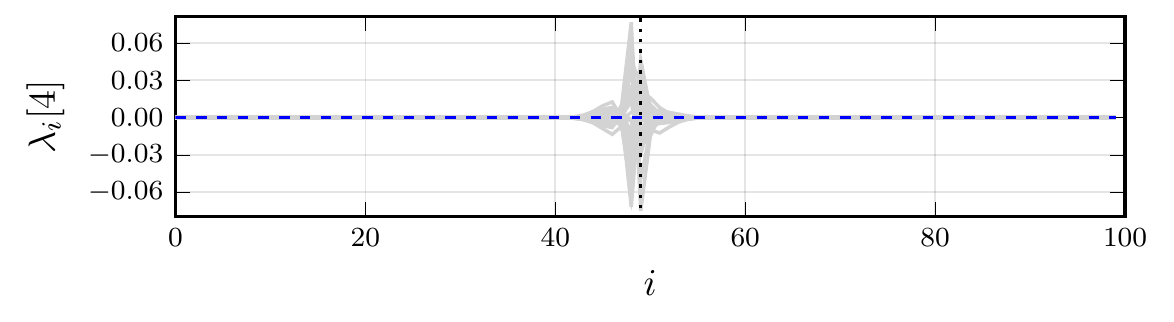}&
    \includegraphics[width=.43\textwidth]{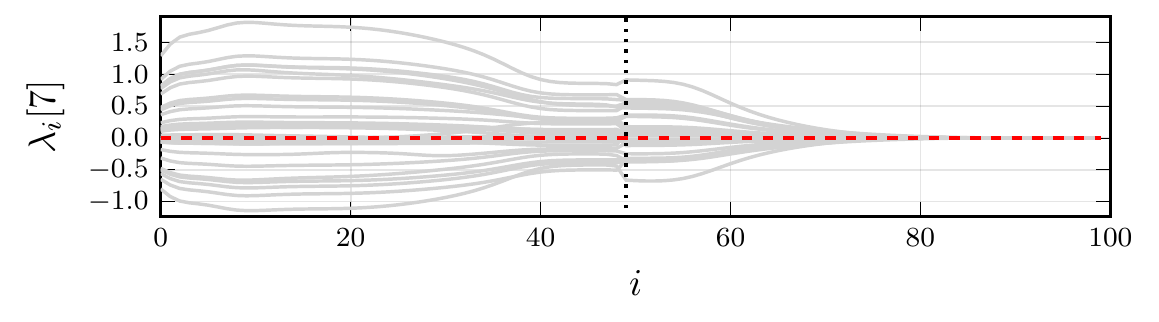}
  \end{tabular}
  \vspace{-0.1in}
  \caption{Base and perturbed solutions. Left: Case 1 ($q=1,b=1$). Right: Case 2 ($q=0,b=0$).}\label{fig:nonlinear}
\end{figure*}

We illustrate the results of  Theorem \ref{thm:eds} and of Corollaries \ref{cor:sys}, \ref{cor:ti}. In this study, we solve the problem with base data $d^\star_{-1:N}$ to obtain the base solution $w^\star_{-1:N}$. We then solve a set of problems with perturbed data; in each of these problems, a random perturbation $\Delta d_j$ is introduced at a selected time stage $j$, while the rest of the data do not have perturbation (i.e., $\Delta d_i=0$ for $i\neq j$). The obtained solutions $w^\dag_{-1:N}(d^\star_{-1:N}+\Delta d_{-1:N})$ for the perturbed problems are visualized along with the base solution $w^\star_{-1:N}$. The scripts can be found here \url{https://github.com/zavalab/JuliaBox/tree/master/SensitivityNMPC}. We consider a quadrotor motion planning problem \citep{hehn2011flying} with the time-invariant setting; the cost functions are given by:
\begin{align*}
  \ell(z;d) :=& (x-d)^\top Q(x-d)+u^\top Ru\\
  \ell_f(x;d) :=& (x-d)^\top Q_f(x-d),\quad\ell_b(x;d)=0,
\end{align*}
where $Q:=\diag(1,1,1,q,q,q,1,1,1)$, $R:=I$, $Q_f:=I$, and $T:=I$; and the dynamic mapping is obtained from:
\begin{subequations}
  \begin{align}
    \frac{d^2X}{dt^2} &= a(\cos\gamma\sin\beta\cos\alpha+\sin\gamma\sin\alpha)\\
    \frac{d^2Y}{dt^2} &= a(\cos\gamma\sin\beta\sin\alpha-\sin\gamma\cos\alpha)\\
    \frac{d^2Z}{dt^2} &= a\cos\gamma\cos\beta-g\\
    \frac{d\gamma}{dt}&= (b\omega_X\cos\gamma + \omega_Y\sin\gamma)/\cos\beta\\
    \frac{d\beta}{dt}&= -b\omega_X\sin\gamma + \omega_Y\cos\gamma\\
    \frac{d\alpha}{dt}&= b\omega_X\cos\gamma\tan\beta + \omega_Y\sin\gamma\tan\beta\ + \omega_Z,
  \end{align}
\end{subequations}
where the state and control variables are defined as: $x:=(X,\dot{X},Y,\dot{Y},Z,\dot{Z},\gamma,\beta,\alpha)$ and $u:=(a,\omega_X,\omega_Y,\omega_Z)$. We use $q$ and $b$ as parameters that influence controllability and observability. In particular, the system becomes less observable if $q$ becomes small and the system loses controllability as $b$ becomes small (the effect of manipulation on $\omega_X$ becomes weak). We have empirically tested the sensitivity behavior for $q=b=1$ (Case 1) and $q=b=0$ (Case 2). One can see that some of the assumptions (e.g., $S_i=0$ in Corollary \ref{cor:sys}) may be violated, but one can also see that, qualitatively, the system is more observable and controllable in Case 1 than in Case 2. The results are presented in Figure \ref{fig:nonlinear}. The base trajectories are shown as dashed lines, the perturbed trajectories are shown as solid gray lines, and the perturbed stages are highlighted using vertical lines. We can see that, for Case 1 ($q=1,b=1$), the differences between the base and perturbed solutions become small as moving away from the perturbation point (EDS holds). On the other hand, for Case 2 ($q=0,b=0$) one cannot observe EDS; this confirms that observability and controllability induce EDS. 

\section{Conclusions}\label{sec:con}
{We have shown that uniform controllability and observability provide sufficient conditions for exponential decay of sensitivity in dynamic optimization. As part of future work, we will aim to establish exponential decay of sensitivity under mesh refinement settings and will aim to establish formal connections with continuous-time results.}

\bibliography{ifacconf}  
\end{document}